\documentclass{article}

\usepackage{arxiv}

\usepackage[utf8]{inputenc} % allow utf-8 input
\usepackage[T1]{fontenc}    % use 8-bit T1 fonts
\usepackage{hyperref}       % hyperlinks
\usepackage{url}            % simple URL typesetting
\usepackage{booktabs}       % professional-quality tables
\usepackage{amsfonts}       % blackboard math symbols
\usepackage{nicefrac}       % compact symbols for 1/2, etc.
\usepackage{microtype}      % microtypography
\usepackage{lipsum}
\usepackage{graphicx}
\graphicspath{ {./images/} }

\title{Approximate Formulas for Characteristics
of Multichannel LIFO Preemptive-Resume Priority Queueing System
}

\author{
 Tatashev A.G. \\
  Department of Higher Mathematics\\
  Moscow Automobile and Road Construction\\
  State Technical University (MADI) \\
  Moscow, Leningradsky avenue, 64, Russia  \\
  \texttt{a-tatashev@yandex.ru} \\
   \And
Seleznjev O.V.\\
Umea University, Sweden\\
  \texttt{oleg.seleznjev@umu.se} \\
   \And
 Yashina M.V. \\
  Department of Higher Mathematics\\
  Moscow Automobile and Road Construction\\
  State Technical University (MADI) \\
  Moscow, Leningradsky avenue, 64, Russia  \\
  \texttt{mv.yashina@madi.ru} \\
}

\begin{document}
\maketitle
\begin{abstract}
This paper considers a multichannel preemptive-resume priority  queueing system with a Poisson input and an arbitrary service time
distribution depending on the priority of job. Jobs of the same priority are serviced according to the LIFO rule. If, at moment of job arrival,
 all servers are busy, and at least one server is busy with the service of a job of a not higher than the priority of arriving job,
 then the service of a job is preempted such that the priority of preempted job is lowest from the priorities of the jobs in service.
 The service of a preempted job is resumed later. The paper proposes approximate formulas for the sojourn time of a prescribed priority job
 and some other characteristics of the system.
\end{abstract}

% keywords can be removed

\section{Introduction}

The papers~\cite{Brosh}--~\cite{Takagi17}
study
a M/M/n queueing system with preemptive-resume priority. It was supposed that the parameter of job service time distribution does not depend on priority, jobs of the same priority are serviced according to the First Input, First Output (FIFO) rule, and the job to be displaced is selected according to Last Come, First Displaced (LCFD) rule. The paper
~\cite{Takagi16} studies an analogous system with Last Input, First Output (LIFO) and First Come, First Displaced (FCFD) rules.

In~\cite{Williams}, ~\cite{Golovina},
 approximate formulas are proposed for characteristics of M/G/m system with non-preemptive priority discipline.
The papers ~\cite{Bondi},~\cite{Maslova} propose different approximate formulas
for the average sojourn time of a prescribed priority job in M/G/m system with preemptive-resume priority discipline.
The paper ~\cite{Maslova}  also proposes approximate formulas for such characteristics of this system as the average waiting time of a prescribed priority job, the probability that the service of a job does not start immediately, average time before the start of the job service provided that this time is not equal to 0, the the average number of job preemptions, the average duration of interruption interval due to a preemption.
In ~\cite{Golovina}--~\cite{Maslova}, it is assumed that the service time distribution depends on the priority of job.
 In ~\cite{Alencar20}, an approximate approach to estimate on-time message delivery probability in a packet-switched network with messages of different priorities.
In ~\cite{Alencar21}, an approximate approach is proposed to evaluate the sojourn time for a queueing system with a limited processor sharing discipline.

This paper proposes approximate formulas for charactristic of M/G/m system with LIFO  preemptive-resume priority discipline with FCFD displacement rule. It is supposed that the service time distribution depends on the priority of job. The paper proposes approximate formulas for the average sojourn time for a job of prescribed priority, the average waiting time, the probability that the service of a job start immediatly, the average time before the start of the job service provided this time is not equal to 0, the average number of job preemptions, the average interruption interval due to preemptions. These formulas are similar to formulas proposed
in ~\cite{Maslova} for an analogous system with FIFO~--- LCFD preemptive-resume priority discipline.

\section{%2
Description of system
}
\label{section:SD}

We use the following notations: $v_i$ is the average sojourn time for a job of the priority $i;$ $w_i$ is the average waiting time including the service interruptions for the priority $i$ job; $p_i$ is the probability that the service of priority $i$ class job does not start immediately; $u_i$ is the average time before the start of the priority $i$ job service provided this time is not equal to 0; $h_i$ is the average number of the priority $i$ job preemptions $i=1,\dots,N;$ $g_i$ is the duration of a service interruption interval for the priority $i$ job, $i=2,\dots,N;$ $\Lambda_{i}$ is the total arrival rate of priority-classes no lower than $i:$ $\Lambda_i=\lambda_1+\dots+\lambda_i;$ $R_i$ is the load due to priority-classes no lower than $i:$ $R_i=(\lambda_1+\dots+\lambda_i)/m,$   $i=1,\dots,N.$

Denote by $c_i$ the probability of non-zero waiting for $M/G/m$ system computed by well-known Erlang's formula for a waiting system with the arrival load $R_i:$
$$c_i=\frac{(mR_i)^m}{m!(1-R_i)\sum\limits_{k=0}^{m-1}
\frac{(mR_i)^k}{k!}+(mR_i)^m}.$$

The  following equalities are true: $$w_1=p_1u_1,\eqno(1)$$
$$w_i=p_iu_i+h_ig_i,\ i=2,\dots,N,\eqno(2)$$
$$v_i=w_i+b_i,\  i=1,\dots,N,\eqno(3)$$
$$h_i=\frac{\Lambda_i(p_i-p_{i-1})}{\lambda_i},\ i=1,\dots,N.\eqno(4)$$
The proof of (4) is similar to the proof of an analogous statement for a preemptive priority system such that, in this system, the service distribution is exponential with average value independent of the priority class.
The proof of (4) is the following. The probability that all servers are busy by jobs of priority-classes not lower $i$ and there is at least one job of priority-class $i$ equals the difference of the probability that all servers are busy by jobs of priorities not lower than $i$ and  the probability that all servers are busy by jobs fo priorities not lower $i-1,$ and therefore this probability is $p_i-p_{i-1}.$ Hence the average number of preemptions of the priority $i$ per a time unit is equal to $\Lambda_{i-1}(p_i-p_{i-1}).$ From this, taking into account that the average number of arriving priority-class $i$ jobs per a time unit is equal to $\lambda_i,$ one gets (1).

\section{%3
Exact formulas for special cases
}
\label{section:SC}

If $m=1$ (a one-channel system), then the following formulas are true
 ~\cite{Jaiswal}:
$$p_i=R_{i-1}=c_{i-1},\ i=1,\dots,N,\eqno(5)$$
$$u_i=\frac{\sum\limits_{j=1}^{i-1}\lambda_jb_j^{(2)}}
{2R_{i-1}(1-R_{i-1})(1-R_i)},\ i=1,\dots,N,\eqno(6)$$
$$h_i=\Lambda_ib_i,\  i=1,\dots,N,\eqno(7)$$
$$g_i=\frac{R_i}{\Lambda_i(1-R_i)},\ i=2,\dots,N,\eqno(8)$$
$$w_i=\frac{\sum\limits_{j=1}^{i-1}\lambda_jb_j^{(2)}}
{2(1-R_{i-1})(1-R_i)}+\frac{R_ib_i}{1-R_i},\ i=1,\dots,N,\eqno(9)$$
$$v_i=\frac{\sum\limits_{j=1}^{i-1}\lambda_jb_j^{(2)}}
{2(1-R_{i-1})(1-R_i)}+\frac{b_i}{1-R_i},\ i=1,\dots,N.\eqno(10)$$

If the service time is distributed exponentionally with the same average value $b,$ then the following formulas  are true
(in ~\cite{Tatashev1984}, similar formulas were get for a FIFO preemptive-resume priority system):
$$p_i=c_{i-1},\ i=1,\dots,N, \eqno(11)$$
$$u_i=\frac{b}{m(1-R_{i-1})(1-R_i)},\  i=1,\dots,N, \eqno(12)$$
$$h_i=\frac{\Lambda_i(c_i-c_{i-1})}{\lambda_i},\ i=1,\dots,N,\eqno(13)$$
$$g_i=\frac{R_i}{\Lambda_i(1-R_i)},\ i=2,\dots,N,\eqno(14)$$
$$w_i=\frac{c_{i-1}b}{m(1-R_{i-1})(1-R_i)}+\frac{R_i(c_i-c_{i-1})}{\lambda_i(1-R_i)},\ i=1,\dots,N, \eqno(15)$$
$$v_i=\frac{c_{i-1}b}{m(1-R_{i-1})(1-R_i)}+\frac{R_i(c_i-c_{i-1})}{\lambda_i(1-R_i)}+b,\ i=1,\dots,N. \eqno(16)$$

\section{%3
Approximate formulas
}
\label{section:AF}

The following approximate formulas for the considered with arbitrary service time distribution may be proposed $(i=1,\dots,N):$
$$p_i=c_{i-1},\  \eqno(17)$$
$$u_i=\frac{\sum\limits_{j=1}^{i-1}\lambda_jb_j^{(2)}}
{2m^2R_{i-1}(1-R_{i-1})(1-R_i)}, \eqno(18)$$
$$g_i=\frac{R_i}{\Lambda_i(1-R_i)},\eqno(19)$$
$$h_i=\frac{\Lambda_i(c_{i}-c_{i-1}}{\lambda_i},\eqno(20)$$
$$w_i=\frac{c_{i-1}\sum\limits_{j=1}^{i-1}\lambda_jb_j^{(2)}}
{2m^2R_{i-1}(1-R_{i-1})(1-R_i)}+\frac{R_i(c_i-c_{i-1})}{\lambda_i(1-R_i)},\eqno(21)$$
$$v _i=\frac{c_{i-1}\sum\limits_{j=1}^{i-1}\lambda_jb_j^{(2)}}
{2m^2R_{i-1}(1-R_{i-1})(1-R_i)}+\frac{R_i(c_i-c_{i-1})}{\lambda_i(1-R_i)}+b_i.\eqno(22)$$

 Let us present arguments in favor of the formulas (17)--(22).

1. The formulas (17)--(22) are exact for $m=1$ (one-channel system).

2.  The formulas (17)--(22) are exact in the case of exponential distribution with average value independent of  priority class.

3. If $m=1$ or the service distribution is exponential with average value independent of priority class, then the probability that the service of a priority $i$ job does not start immediately, is equal to the probability that all servers of a related non-preemptive system with arriving load  $R_{i-1}$
Assume this relation holds approximately for the studied multichannel priority system and the related non-priority M/G/m system.
As it noted in ~\cite{Jaiswal}, the probability that all servers of a M/G/m system not sufficient depends on the job service distribution provided that the averaged service time is presribed. Therefore this probability is equal to the probability for the related system M/M/m approximately, i.e., we get (14).  Note that the analogous formulas are proposed
in~\cite{Williams} and ~\cite{Golovina} for  related systems with non-preemptive priorities.

4. If $m=1$ or the service time distribution is exponential with average value indepenent of priority class, then the average time interval before the priority class $N$ the job service start is the same for the considered system and the related non-priority system with priorities classes $i=1,\dots,N-1.$ Assume that the average time intervals before the priority class $N$ the job service start are approximately equal to each other for the studied multichannel preemptive-resume priority system and the elated non-priority system in the case of multichannel systems. Suppose $d_N=p_Nu_N,$ where the values $p_N$ and $u_N$ are computed approximately according to (14), (15). Then the value of $d_i$ is the same as the value computed according to an approximate formula for the related non-priority discipline proposed in ~\cite{Hokstad}. This is an argument in favor of formula (15).

5. If $m=1$ or the job service time distribution is exponential with average value independent of priority-class $(i\ge 2),$ then the average interruption interval $g_i$ will be the same for the considered system and for the related non-priority system with the service rate value multiplied by $m.$ Assume that this relation holds approximately for the considered multichanel priority system, we get the equation (7).

6. Combining (1)--(4) and (17)--(19), we get (20)--(22).

7. If the job service time is exponential with parameter that may depend on priority class, then the values of sojourn time are the same for the FIFO preemptive-resume discipline and
 the LIFO preemptive-resume discipline. In accordance with this, the formulas (21), (22) are the same as the formulas for the waiting time
and sojourn time in the related system with the FIFO preemptive-resume discipline.
 \vskip 20pt
Assume that the job service time priority is exponential with average value depending on priority-class:
$$B_i(x)=1-e^{-\mu_ix},\ i=1,2,3,4,$$
$$\mu_1=5,\ \mu_2=\frac{5}{2},\ \mu_3=\frac{5}{3},\ \mu_4=\frac{5}{4},$$
$$m=3,\ \lambda_1=\lambda_2=\lambda_3=\lambda_4=1.$$
Then the value computed according to formulas (21), (22) are
$$w_1=0.000084,\ w_2=0.0053,\ w_3=0.075,\ w_4=0.65,$$
$$v_1=0,20,\ w_2=0.40,\ w_3=0.61,\ w_4=1.4.$$
The values obtained by simulation are
$$w_1=0,000084,\ w_2=0.0054,\ w_3=0.075,\ w_4=0.60,$$
$$v_1=0.20,\ w_2=0.40,\ w_3=0.61,\ w_4=1.4.$$

\bibliographystyle{unsrt}


\begin{thebibliography}{1}


\bibitem{Brosh}
%1.
 Brosh I. Preemptive priority assignment in multichannel systems. Operations Research, 1969, vol. 17, pp.~526–535.

\bibitem{Segal}
%2.
 M. Segal, A multiserver system with preemptive priorities, Opns. Res., 18 (1970), 316–323,
\newblock{ DOI: 10.1287/opre.18.2.316}

\bibitem{ Buzen}
%3.
 Buzen J.P., Bondi A.B. The response times of priority classes under preemptive resume
in M/M/m queues, Operations Research, 31 (1983), 456–465.


\bibitem{Tatashev1984}
%4.
Tatashev A.G. Calculation of the distribution of the waiting time in a multiple-channel
queueing system with fixed priorities, Engineering Cybernetics, 22 (1984), 59–62.
 (Originally published in Tekhnicheskaya Kibernetika),
 6 (1983), 163–166.


\bibitem{Zeltyn }
%5.
  Zeltyn S., Feldman Z.,  Wasserkrug S. Waiting and sojourn times in a
multi-server queue with mixed priorities, Queueing Systems, 61 (2009), 305–328,
\newblock {DOI:10.1007/s11134-009-9110-4}

\bibitem{Takagi17}
%6.
Takagi H. Unified and refined analysis of the response time and waiting time in the M/M/m FCFC preemptive-resume priority queue. Journal of Industrial and Management Optimization, vol.~13, no.~4, 2017, pp.~1945--1973.
\newblock { DOI:10.3934/jimo.2017026}

\bibitem{Takagi16}
% 7.
 Takagi H. Analysis of the response and waiting times in the M/M/m LCFC     preemptive-resume priority queue. International
Journal of Pure and Applied Mathematics, vol. 109, no. 2, 2016, pp. 325--370.
DOI: 10.12732/ijpam.v109i2.12

\bibitem{Williams}
%8.
Williams T. Nonpreemptive multy-server-priority queues. J. of the Operat. Res. Soc., 1980, vol.~31, no.~12, pp.~1105--1108.

\bibitem{Golovina}
%9.
Golovina E.V., Dukhovnyi I.M. On efficacy of approximate approach to evaluation of complex queueing systems based on use of relations invariants. In: Information systems and automatic commutation. Moscow, 1975. P. 13--14.


\bibitem{Bondi}
%10.
 Bondi A.P., Busen J.P. The response times of priority classes under preemptive
resume in M/G/m queues. Peform. Eval. Rev., 1984, vol. 12, no. 3, pp. 195–201.

\bibitem{Maslova}
%11.
 Maslova E.A., Tatashev A.G. Approximate method for preemptive priority service
in a multichannel system. Automatic Control and Computer Science, 1992, vol. 26, no 1,
pp. 10–15. Translated from Avtomatika i Vychislitelnaya Tekhnika, 1992, no. 1, pp. 13–18.

\bibitem{Alencar20}
%12.
 Alencar M.S., Tatashev A.G., Yashina M.V.. Stochastic estimation of on-time message delivery in a packet-switched network with messages of different priorities.  2020 International Conference on Engineering Management of Communication and Technology (EMCTECH), 2020, pp. 1-4,\\
\newblock {DOI: 10.1109/EMCTECH49634.2020.9261511.}

\bibitem{Alencar21}
%13.
 Alencar M.S., Buslaev D.A., Tatashev A.G., Yashina M.V. Estimation method of  v2v networks  capacity  in saturated traffic  flows.
2021 Systems of signals generating and processing in the field of on board communications.
IEEE Conference \#51389. 16--18 March 2021, Moscow, Russia.

\bibitem{Jaiswal}
%14.
 Jaiswal N.K. Priority queues. Vol. 50 of Mathematics and Science Engineering. Academic Press, New York and London, 1968.

\bibitem{Stoyan}
%15.
 Stoyan D. Qualitative Eigenschaften und Abschtzungen stochastischer Modelle. Published by Gruyter.  Berlin, 1977.
 \newblock {DOI: 10.1515/9783112563762}

\bibitem{Hokstad}
%16.
 Hokstad P. Approximations for M/G/m queue. Oper. Res., 1978, vol. 10, no.~1, pp.~510--523.




\end{thebibliography}

\end{document}